\documentclass[11pt]{amsart}
\newtheorem{theorem}{Theorem}[section]

\newtheorem{proposition}[theorem]{Proposition}

\newtheorem{lemma}[theorem]{Lemma}

\def\qed{\hfill $\Box$\medskip}

\def\IC{{\bf C}}

\def\cS{{\mathcal S}}

\def\tr{{\rm tr}\,}

\def\span{{\rm span}\,}
\def\conv{{\rm conv}\,}
\def\cE{{\mathcal E}}

\begin{document}
\openup .9\jot
\title{Numerical Ranges of the product of Operators}
\author{Hongke Du, Chi-Kwong Li, Kuo-Zhong Wang, Yueqing Wang, Ning Zuo}

\date{}
\maketitle

\begin{abstract}
We study containment regions
of the numerical range of the product of operators $A$ and $B$ such that
$W(A)$ and $W(B)$ are line segments. It is shown that the containment region is equal to
the convex hull of elliptical disks determined by the spectrum of $AB$,
and conditions on $A$ and $B$ for the set equality holding are obtained.
The results cover the case when $A$ and $B$ are self-adjoint operators  extending
the previous results on the numerical range of the product of two orthogonal projections.
\end{abstract}

Keywords: Numerical range, product of matrices and operators.

AMS Classification. 15A60, 47A12.

\section{Introduction}

Let $B(H)$ be the algebra of bounded linear operators on a complex Hilbert space $H$.
We identify $B(H)$ with $M_n$, the algebra of $n$-by-$n$ complex matrices, if $H$ has
finite dimension $n$.
The spectrum $\sigma (A)$, and the numerical range $W(A)$ of an operator
$A\in B(H)$ are defined by
$$\sigma(A) =
\{\lambda : A-\lambda I\mbox{ is not invertible}\}
\quad \mbox{ and } \quad W(A)=\{\langle Ax,x,\rangle :x\in H\; , \|x\|=1\},
$$
respectively.
Here $\langle \cdot ,\cdot\rangle$ and $\|\cdot\|$ are the standard inner product
and its associated norm on $H$, respectively.
The spectrum and the numerical range are useful tools in the study of matrices
and operators; for example, see \cite{K,RA,PR}.
It is known that $W(A)$ is a bounded convex subset of $\IC$.
When $H$ is finite dimensional, it is compact.
In general, the closure of the numerical range satisfies
 $\sigma (A)\subseteq \overline{W(A)}$.
Especially, for $A\in M_2$, $W(A)$ is an elliptical disk
with $\lambda_1$ and $\lambda_2$ as foci and
$\{\mbox{tr}(A^*A)-|\lambda_1|^2-|\lambda_2|^2\}^{1/2}$ as
minor axis, where $\lambda_1$ and $\lambda_2$ are eigenvalues of $A$.

An operator $A\in B(H)$ is an orthogonal projection if $A^2=A=A^*$, contraction if
$\|A\|\equiv \sup_{\|x\|=1}\|Ax\|\le 1$, and positive if $\langle Ax,x\rangle\ge 0 $ for
all $x\in H$.
In \cite{Kl}, it was shown that if $P,Q\in B(H)$ are orthogonal projections and
$0\in \sigma (P)\cup \sigma (Q)$, then
$$
\overline{W(PQ)}=\overline{\mbox{conv}\{\cup_{\lambda\in \sigma (PQ)}\cE(\lambda)\}},
$$
where $\mbox{conv}\{\cS\}$ is the convex hull of the set $\cS$
and $\cE(\lambda)$ is the ellipse disc with foci $0$ and $\lambda$, and
length of minor axis $\sqrt{\lambda (1-\lambda)}$.
In general, the following example shows that
the above equality may not hold for positive contractions
$A,B\in B(H)$.
Let $A=B={\small \left(
           \begin{array}{cc}
             1 & 0 \\
             0 & 1/2 \\
           \end{array}
         \right)}$.
Then $W(AB)=[1/4,1]\neq \overline{\mbox{conv}\{\cE(1)\cup \cE(1/4)\}}$.

In this paper, we consider the containment regions for the
numerical range of the product of positive contractions, and extend the result
to more general operators, namely, those operators with numerical
ranges equal to line segments.

First of all, for two positive contraction operators $A,B\in B(H)$, it is known
that
$$\overline{W(AB)} \subseteq \{x+iy: -1/8 \le x \le 1,  -1/4 \le y \le 1/4\};$$
see \cite{Kl}, \cite{W}, \cite{O}. To see this\footnote{Li would like to thank
Professor Fuzheng Zhang
for showing him this proof.}, we use the positive definite ordering that
$X\ge Y$, and show that

\medskip\centerline{
(a) \ $-I/8 \le (AB+BA)/2 \le I$ \quad and \quad (b)
$-I/4 \le (AB-BA)/(2i)\le I/4$.}

\medskip\noindent
For (a), it is clear that $\|AB+BA\| \le 2$ so that $-2I \le AB+BA \le 2I$
for two positive contractions $A$ and $B$. Furthermore, note that
$A+B - A^2 - B^2 \ge 0$, and hence
\begin{eqnarray*}
0 &\le& (A+B-I/2)^2  = A^2 + B^2 + (AB+BA) + I/4 - A - B\\
&=& (AB+BA) + I/4 + (A^2-A) + (B^2-B) \le (AB+BA) + I/4.\end{eqnarray*}
For (b), since $\|A-I/2\|\le 1/2$ and $\|B-I/2\|\le 1/2$, we have
\begin{eqnarray*}
\|i(AB-BA)\|&=&\|(A-I/2)(B-I/2)-(B-I/2)(A-I/2)\|\\
&\le &2\|(A-I/2)(B-I/2)\|\le 1/2.\end{eqnarray*}
Moreover, $i(AB-BA)$ is self-adjoint and then condition $(b)$ holds.

\medskip
Suppose $\lambda \in [0,1]$. Denote by $\cE(\lambda)$
the elliptical disk with foci $0,\lambda$, minor axis with end points
$(\lambda \pm  i\sqrt{\lambda(1-\lambda)})/2$, and major
axis with end points $(\lambda \pm \sqrt{\lambda})/2$.
Then $W\left({\tiny \left(
          \begin{array}{cc}
            \lambda & 0 \\
            \sqrt{\lambda(1-\lambda)} & 0 \\
          \end{array}
        \right)}
\right) = \cE(\lambda).$
We have the following result in \cite{K}.

\begin{theorem} \label{1.1}
Let $P, Q \in M_n$ be non-scalar orthogonal projections. Then
$$W(PQ) = \conv\{\cup_{\lambda \in \sigma(PQ)} \cE(\lambda)\}.$$
\end{theorem}

One can obtain the above result
using the following canonical  for a product of projections $P, Q \in M_n$;
see \cite{Kl,S} and its references.

\begin{proposition} \label{1.2}
Suppose $P, Q \in M_n$ are non-scalar projections and $U \in M_n$ is
unitary such that
$U^*(P+iQ)U$ is a direct sum of $(I_p + iI_p) \oplus I_{q} \oplus iI_{r} \oplus 0_{s}$,
and
$$C_j = \begin{pmatrix}c_j^2 + i & c_js_j \cr c_js_j & s_j^2\cr\end{pmatrix}, \quad j = 1, \dots, k,$$
where $c_j \in (0,1), s_j = \sqrt{1-c_j^2}$.
Then $U^*PQU$ will be a direct sum of $I_p \oplus 0_{q+r+s}$ and
$$\hat C_j = \begin{pmatrix}c_j^2 & 0 \cr c_js_j & 0 \cr\end{pmatrix},
\quad j = 1, \dots, k.$$
\end{proposition}

In the next two sections, we will consider containment regions 
for the numerical range of the product of 
a pair of positive contractions, and extend the results to a more general class of matrices,
namely, those matrices with numerical ranges contained in line segments.
We will consider the infinite dimensional version of Theorem \ref{1.1} and its
generalization in Section 4.

\section{Positive Contractions}

In this section, we extend Theorem \ref{1.1} to obtain a containment region $\cS$ of $W(AB)$
for two positive contractions
$A, B \in M_n$, and determine the conditions for $\cS = W(AB)$.
We begin with some technical lemmas. We will denote by
$\lambda_1(X) \ge \cdots \ge \lambda_n(X)$ the eigenvalues of a Hermitian matrix
$X\in M_n$.

\begin{lemma} \label{2.0}  Suppose $T = T_1 \oplus \dots \oplus T_m
\oplus T_0 \in M_n$ and $\mu = \mu_1 + i \mu_2 \in \IC$ satisfying
$T_1 = \cdots = T_m  \in M_2$ are non-scalar matrices,
$$2\mu_1 = \lambda_1(T_1+T_1^*) > \lambda_1(T_0).$$
Then up to a unit multiple, there is a unique unit vector $\tilde{x} \in \IC^2$ such that
$\tilde{x}^* T_1 \tilde{x} = \mu_1 + i \mu_2$. Moreover,
if $x\in \IC^n$ such that  $x^*Tx = \mu_1 + i \mu_2$, then
$$x = v\otimes \tilde x \oplus 0_{n-2m} =
(v_1\tilde x^t, \dots, v_m \tilde x^t, 0, \dots, 0)^t$$
where $v = (v_1, \dots, v_m)^t\in \IC^m$ is a unit vector.
\end{lemma}

\it Proof. \rm
Obviously, $\mu\in \partial W(T_1)$.
Hence there is a unit vector $\tilde{x}\in \IC^2$ such that
$\tilde{x}^* T_1 \tilde{x} = \mu_1 + i \mu_2$.

Now, suppose $x = x_1 \oplus \cdots \oplus x_m \oplus x_0$, where,
$x_1, \dots, x_m \in \IC^2$ and $x_0 \in \IC^{n-2m}$,
is a unit vector such that $x^*Tx = \mu_1 + i \mu_2$.
Then
$$2\mu_1 = x^*(T+T^*)x =  \sum_{j=0}^m x_j^*(T_j+T_j^*)x_j  \le
2\mu_1 \sum_{j=1}^m \|x_j\|^2 + \lambda_1(T_0+T_0^*) \|x_0\|^2.$$
Thus, $x_0 = 0$ and $(T_j+T_j^*)x_j = 2\mu_1 x_j$ for $j = 1, \dots, m$.
and $x_j = v_j \tilde x$ with $v_j \in \IC$ for $j = 1, \dots, m$.
Let $v = (v_1, \dots, v_m)^t$. Then $x = v \otimes \tilde x \oplus 0_{n-2m}$,
and $\|v\| = \|x\|/\|\tilde x\| = 1$ as asserted.
\qed

\begin{lemma} \label{2.1}
Let $P, Q \in M_n$ be non-scalar orthogonal projections.
Suppose
that there is a supporting line $L$ of $W(PQ)$ satisfying
$L \cap W(PQ) = \{\mu\} \subseteq \cE(\hat \lambda)$ with $\hat \lambda \in \sigma(PQ)$ and $\hat{\lambda}\in (0,1)$.
Suppose that $\mu\notin \cE (\lambda)$ for all other $\lambda\in \sigma (PQ)$ and that $\langle PQx,x\rangle=\mu$ for some unit vector $x$.
Then there is a unitary matrix $V \in M_n$ with the first two columns $v_1, v_2$
such that $\span\{v_1,v_2\}=\span\{x,PQx\}$, and
$$V^*PV ={\small
\left(\begin{array}{cc} \hat \lambda & \sqrt{\hat\lambda - \hat \lambda^2} \\
\sqrt{\hat\lambda - \hat \lambda^2} & 1-\hat\lambda \\ \end{array}  \right)} \oplus P'
\qquad   \hbox{ and } \qquad V^*QV =
{\small \left(\begin{array}{cc} 1 & 0 \\ 0 & 0 \\ \end{array}  \right)} \oplus Q'.$$
\end{lemma}

\it Proof. \rm
Suppose $PQ$ has the canonical form described in Proposition \ref{1.2}, and
$U$ is the unitary such that $U^*PQU = C_1 \oplus \cdots \oplus C_k \oplus I_r \oplus 0_s$,
$P=P_1\oplus \cdots \oplus P_k \oplus I_r \oplus P' $
and $Q=Q_1\oplus \cdots \oplus Q_k \oplus I_r \oplus Q'$, where $P'Q'=0_s$ and
$C_i=P_iQ_i$ for $1\le i\le k$.
We may further assume that $C_1, \dots, C_m$ satisfy
$W(C_1) = \cdots = W(C_m) = \cE(\hat \lambda)$ such that
$\mu \notin W(C_j)$ for all other $j\in \{m+1,\ldots,k\}$.
Thus,

\medskip
$C_1=\cdots=C_m=\begin{pmatrix}c^2 & 0 \cr cs & 0 \cr\end{pmatrix}$,
$P_1=\cdots=P_m =
\left(\begin{array}{cc} c^2 & cs \\ cs & s^2 \\ \end{array}  \right)$ and
$Q_1=\cdots=Q_m=\left(\begin{array}{cc} 1 & 0 \\ 0 & 0 \\ \end{array} \right)$

\medskip

\noindent
with $\hat{\lambda} = c^2$ and $s^2 = \sqrt{1-c^2}$.
Because $L \cap W(PQ) = \{\mu\}$, there is
$t \in [0, 2\pi)$ such that $e^{it}\mu + e^{-it}\bar \mu$  is
the largest eigenvalue of $e^{it}PQ + e^{-it}QP$;
e.g., see \cite[Chapter 1]{RA}.
Now, set $\hat P + i\hat Q = U^*(P+iQ)U$ and $\hat x = U^*x$ so that
$\hat x^* \hat P \hat Q \hat x = \mu$. Then  $T = e^{it}\hat{P}\hat{Q}$
will satisfy the hypothesis of Lemma \ref{2.0}.
It follows that
$$\hat x = U^*x = v_0 \otimes \tilde x \oplus 0_{n-2m} \quad \hbox{ and  } \quad
\hat P \hat Q \hat x = (U^*PQU)(U^*x) = v_0 \otimes \tilde y \oplus 0_{n-2m}\equiv \hat{y},$$
where $v_0\in \IC^m$ with $\|v_0\|=1$ and $\tilde y = C_1 \tilde x$.
For any $y\in \mbox{span}\{\hat{x},\hat{y}\}$, we have $y=v_0\otimes y_0\oplus 0_{n-2m}$
for some $y_0\in \mbox{span}\{\tilde{x},\tilde{y}\}$, and then
$$\hat{P}\hat{Q}y=v_0\otimes C_1y_0\oplus 0_{n-2m},\;
\hat{P}y=v_0\otimes P_1y_0\oplus 0_{n-2m},\; \hat{Q}y=v_0\otimes Q_1y_0\oplus 0_{n-2m}.\leqno{(1)}
$$
Let $V'\in M_n$ be a unitary such that the span of the first two columns of $V'$
contains the set $\{x,PQx\}$, and let
$\hat{V}=(\hat{v}_1,\ldots,\hat{v}_n)=U^*V'$.
Then $\hat{V}$ is a unitary and
$\mbox{span}\{\hat{v}_1,\hat{v}_2\}=\mbox{span}\{\hat{x},\hat{y}\}$.
From $(1)$, we obtain that
$$\hat{V}^*\hat{P}\hat{Q}\hat{V}=C'_1\oplus C',\;\hat{V}^*\hat{P}\hat{V}=P'_1\oplus P',
\quad \mbox{ and } \quad
\hat{V}^*\hat{Q}\hat{V}=Q'_1\oplus Q',$$
where $C_1'\cong C_1
$ and $P_1',Q_1'$ are two 2-by-2 orthogonal projections.
Since $c^2=\hat{\lambda}\in (0,1)$ and $P_1'Q_1'=C_1'$,
$Q_1'\neq 0_2$, $I_2$.
There is a unitary $\hat{R}\in M_2$ such that
$$
\hat{R}^*Q_1'\hat{R}=\left(
                       \begin{array}{cc}
                         1 & 0 \\
                         0 & 0 \\
                       \end{array}
                     \right)
                     \quad \mbox{ and } \quad
\hat{R}^*P_1'\hat{R}=\left(
                       \begin{array}{cc}
                         p_{11} & p_{12} \\
                         \overline{p}_{12} & p_{22} \\
                       \end{array}
                     \right).
$$
Hence $\left(
         \begin{array}{cc}
           p_{11} & 0 \\
           \overline{p}_{12} & 0 \\
         \end{array}
       \right)=\hat{R}^*C_1'\hat{R}\cong C_1
$,
and then $c^2=p_{11}, \overline{p}_{12}=e^{i\theta}cs$ for some $\theta\in [0,2\pi)$.
Let $R=\left(\hat{R}\left(
                   \begin{array}{cc}
                     1 & 0 \\
                     0 & e^{i\theta} \\
                   \end{array}
                 \right)
\right)\oplus I_{n-2}$, and $V = U\hat{V}R$.
Then
$$V^*PV  = R^*\hat V^*\hat P \hat VR =  {\small
\left(\begin{array}{cc} \hat \lambda & \sqrt{\hat\lambda - \hat \lambda^2} \\
\sqrt{\hat\lambda - \hat \lambda^2} & 1-\hat\lambda \\ \end{array}  \right)} \oplus P'
\ \hbox{ and } \
V^*QV = R^*\hat V^* \hat Q \hat VR =
{\small \left(\begin{array}{cc} 1 & 0 \\ 0 & 0 \\ \end{array}  \right)} \oplus Q'$$
as asserted.
\qed

\begin{theorem} \label{2.2}
Let $A,B \in M_n$ be two non-scalar positive contractions. Then
$$W(AB) \subseteq \conv\{\cup_{\lambda \in \sigma(AB)} \cE(\lambda)\}.$$
The set equality holds if and only if
there is a unitary matrix $U$ such that
$U^*AU  = A' \oplus A^{''}, U^*BU = B' \oplus B^{''}$ such that
$A', B'$ are orthogonal projections such that
$W(A^{''}B^{''})\subseteq W(A'B') = W(AB)$.
\end{theorem}

\it Proof. \rm
Let
$$\hat A = \begin{bmatrix}
A & \sqrt{A-A^2} & 0 \cr
\sqrt{A-A^2} & I_n-A& 0\cr
0 & 0 & 0 \cr
\end{bmatrix} \quad \hbox{ and } \quad
\hat B = \begin{bmatrix}
B & 0 & \sqrt{B-B^2} \cr
0 & 0 & 0 \cr
\sqrt{B-B^2} & 0 & I_n-B\cr
\end{bmatrix}.$$
Then
$$T = \hat A\hat B =
\begin{bmatrix}
AB & 0 & A\sqrt{B-B^2} \cr
\sqrt{A-A^2}B & 0 & \sqrt{(A-A^2)(B-B^2)}\cr
0 & 0 & 0 \cr
\end{bmatrix}$$
satisfies
$\sigma(\hat A\hat B) = \sigma(AB) \cup \{0\}$ and
$$W(AB) \subseteq W(T) = \conv\{\cup_{\lambda \in \sigma(\hat A\hat B)}\cE(\lambda)\}.$$

Now, suppose that $W(AB)=\conv\{\cup_{\lambda \in \sigma(AB)} \cE(\lambda)\}=W(T)$.
Then
$$W(AB)=\conv \{S \cup_{\lambda\in \sigma (AB)\setminus S}\cE(\lambda)\},$$
where $S=\sigma (AB)\cap \{0,1\}$.
Obviously, $\sigma (AB)\setminus S=\emptyset$ if and only if $W(AB)\subseteq [0,1]$.

If $W(AB)\subseteq [0,1]$, then $AB$ is a Hermitian matrix so that $AB=B^*A^*=BA$.
Hence $A$ and $B$ commute, and there is a unitary
$U$ such that $A=U^*\Lambda_1U$ and $B=U^*\Lambda_2U$,
where $\Lambda_1=\mbox{diag}(a_1,\ldots,a_n)$ and $\Lambda_2=\mbox{diag}(b_1,\ldots,b_n)$.
Then $W(AB)=W(\Lambda_1\Lambda_2)=[\alpha_0,\alpha_1]$, where $\alpha_0=\min_{1\le i\le n}a_ib_i$ and $\alpha_1=\max_{1\le i\le n}a_ib_i$.
Hence we have the desired conclusion.

Next, suppose that $W(AB)$ is not in $ [0,1]$.
This is $\sigma (AB)\setminus S\neq \emptyset$.
Let $\lambda_1\in \sigma (AB)\setminus S$ be such that $\partial \cE(\lambda_1)\cap \partial W(AB)$ contains an arc.
Then there exists $\mu\in \partial \cE(\lambda_1)\cap \partial W(AB)$ with $\mu\notin \cE(\lambda)$ for all other $\lambda\in \sigma (AB)$.
Let $x_1\in \IC^n$ be a unit vector with $x_1^*ABx_1=\mu$.
Since $\partial W(AB)=\partial W(T)$, there is $\theta_1\in [0,2\pi)$ satisfying
$2\mbox{Re}(e^{i\theta_1}\mu)=\max \sigma (e^{i\theta_1}T+e^{-i\theta_1}T^*)$.
Let $\hat{x}_1=x_1\oplus 0_{2n}$.
Then $\hat{x}_1$ is an eigenvector of $e^{i\theta_1}T+e^{-i\theta_1}T^*$ corresponding
to $e^{i\theta_1}T$ so that
$$(e^{i\theta_1}AB+e^{-i\theta_1}BA)x_1
=2\mbox{Re}(e^{i\theta_1}\mu)x_1,\;e^{i\theta_1}\sqrt{A-A^2}Bx_1=0,\mbox{ and }
e^{-i\theta_1}\sqrt{B-B^2}Ax_1=0.$$
Hence $T\hat{x}_1=ABx_1\oplus 0_{2n}$.
By Lemma 2.2, there is a unitary $\hat{U}_1$ and $\hat{U}_1=U_1\oplus I_{2n} $ such that
the span of the first two columns of $U_1$ contains the set $\{x_1,ABx_1\}$,
$$
\hat{U}_1^*\hat{A}\hat{U}_1=\left(
                              \begin{array}{cc}
                                \lambda_1 & \sqrt{\lambda_1(1-\lambda_1)} \\
                                \sqrt{\lambda_1(1-\lambda_1)} & 1-\lambda_1 \\
                              \end{array}
                            \right)\oplus  \hat{A}_1
\quad \mbox{ and } \quad
\hat{U}_1^*\hat{B}\hat{U}_1=\left(
                              \begin{array}{cc}
                                1 & 0 \\
                                0 & 0 \\
                              \end{array}
                            \right)\oplus \hat{B}_1,
$$
where
$$
\hat{A}_1=\left[
                                            \begin{array}{ccc}
                                              A'_1 & C_1^* & 0 \\
                                              C_1 & I-A & 0 \\
                                              0 & 0 & 0 \\
                                            \end{array}
                                          \right]\quad
\mbox{ and } \quad \hat{B}_1=\left[
                                            \begin{array}{ccc}
                                              B'_1 & 0 & D_1^* \\
                                              0 & 0 & 0 \\
                                              D_1 & 0 & I-B \\
                                            \end{array}
                                          \right].$$
Thus,

\medskip
\centerline{$U_1^*AU_1=\left(
                              \begin{array}{cc}
                                \lambda_1 & \sqrt{\lambda_1(1-\lambda_1)} \\
                                \sqrt{\lambda_1(1-\lambda_1)} & 1-\lambda_1 \\
                              \end{array}
                            \right)\oplus A'_1$, \quad
                            $U_1^*BU_1=\left(
                              \begin{array}{cc}
                                1 & 0 \\
                                0 & 0 \\
                              \end{array}
                            \right)\oplus B'_1$,}
\medskip\noindent
and
$U_1^*ABU_1=C_1\oplus A'_1B'_1
$,
where $C_1={\small \left(
                \begin{array}{cc}
                  \lambda_1 & 0 \\
                  \sqrt{\lambda_1(1-\lambda_1)} & 0 \\
                \end{array}
              \right)}.
$
Then $A'_1,B'_1$ are positive contractions, $\hat{A}_1,\hat{B}_1$ are orthogonal projections,
and
$$\conv \{ \cup_{\lambda\in \sigma (AB)\setminus S}\cE(\lambda)\}=\conv \{\cE (\lambda_1)\cup_{\lambda\in \sigma (A_1B_1)\setminus S}\cE(\lambda)\} .$$

Now, suppose $W(AB) = W(T) = \conv\{S \cup_{j=1}^k W(C_j)\}$
for $k$ distinct matrices $C_1, \dots, C_k \in M_2$
such that for $j = 1, \dots, k$, $\lambda_j\in \sigma (AB)\setminus S$,
$\partial \cE(\lambda_j)\cap \partial W(AB)$ contains an arc, and $W(C_j) = \cE(\lambda_j)$.
Since the argument in the preceding paragraph is true for any $\cE(\lambda_j)$
for $j = 1,\dots, k$,
there is an orthonormal set $\{v_1, \dots, v_{2k}\} \subseteq \IC^{n}$
and a unitary $V  = V_1 \oplus V_2 \in M_{3n}$, where the first $2k$ columns
of $V_1$ equals $v_1, \dots, v_{2k}$,
such that $V^*TV = C_1 \oplus \cdots \oplus C_k \oplus T_0$.
Thus, $V^*\hat A V = A_1 \oplus \cdots \oplus A_k \oplus \tilde A_0$,
$V^*\hat B V = B_1 \oplus \cdots \oplus B_k \oplus \tilde B_0$.
Consequently, $V_1^*AV_1 = A_1 \oplus \cdots \oplus A_k \oplus A_0$,
and $V^*B V = B_1 \oplus \cdots \oplus B_k \oplus B_0$.

Evidently,  $0 \in \sigma(A_1B_1)\cap \{0,1\} \subseteq S$.
If $1 \notin S$, then $W(AB) = W(A_1B_1)$ so that the conclusion of the theorem
holds with $(A',B') = (A_1,B_1)$.
Suppose that $1 \in S$. Then $1 \in \sigma(A_0B_0)$ because
$\sigma (C_j)=\{0,\lambda_j\}$ with $\lambda_j\in (0,1)$.
Because $A_0,B_0$ are positive contractions,
there is a unitary $U_0$ satisfying
$U_0^*A_0U_0=[1]\oplus A_0'$ and $U_0^*B_0U_0=[1]\oplus B_0'$.
Let $U=(V_1(I_{2k}\oplus U_0))\oplus V_2$. Then
$(U^*AU,U^*BU)  = (A' \oplus A'', B' \oplus A'')$
with $(A',B') = (A_1 \oplus[1], A_2\oplus [1])$,
and the desired conclusion follows.
\qed

\section{Essentially Hermitian matrices}

Recall that a matrix $A \in M_n$ is an essentially Hermitian matrix
if $e^{it}(A - (\tr A)I_n/n)$ is Hermitian for some $t \in [0, 2\pi)$.
\iffalse
\footnote{In operator theory, one may define $A \in B(H)$ to be essentially self-adjoint
if $A - T$ is self-adjoint for a compact operator $T$.  Our definition here is different.}
\fi

\medskip
It is known and not hard to show that the following conditions are equivalent for
$A\in M_n$.

\medskip
\quad (a)  $A$ is essentially Hermitian

\medskip
\quad (b) $W(A)$ is a line segment in $\IC$ joining two complex numbers $a_1, a_2$.

\medskip
\quad (c) $A$ is normal and all its eigenvalues lie on a straight line.

\medskip
The results in the previous section can be extended to essentially Hermitian matrices.
We begin with the following result which follows readily from Proposition \ref{1.2}.

\begin{proposition} \label{3.1}
Suppose $A, B\in M_n$ are normal matrices with
$\sigma(A) = \{a_1, a_2\}$ and $\sigma(B) = \{b_1, b_2\}$.
Then $A = (a_1 - a_2)P + a_2 I_n$ and $B = (b_1 - b_2)Q + b_2 I_n$,
where $P$ and $Q$ are orthogonal projections, and there is a unitary matrix
$U$ such that
$U^*(P+iQ)U$ is a direct sum of $(I_p + iI_p) \oplus I_{q} \oplus iI_{r} \oplus 0_{s}$, and
$$\begin{pmatrix}c_j^2 + i & c_js_j \cr c_js_j & s_j^2\cr\end{pmatrix}, \quad
j = 1, \dots, k,$$
where $c_j \in (0,1), s_j = \sqrt{1-c_j^2}$.
Consequently, $U^*ABU$ is a direct sum of a diagonal matrix $D$ with
$\sigma(D) \subseteq \{a_1b_1, a_1b_2, a_2b_1, a_2b_2\}$ and
$$C_j = \begin{pmatrix}
a_1c_j^2+a_2s_j^2 & (a_1-a_2)c_js_j \cr (a_1-a_2)c_js_j & a_1 s_j^2 + a_2c_j^2\cr\end{pmatrix}
\begin{pmatrix} b_1 & 0 \cr 0 & b_2\cr\end{pmatrix},
\quad j = 1, \dots, k.$$
\end{proposition}

Suppose $A,B \in M_n$ satisfy the hypotheses of the above proposition.
Then
$$W(AB) = \conv\{\cup_{j=1}^k W(C_j) \cup W(D)\}.$$
Evidently, $W(D)$ is the convex hull of the diagonal entries of $D$.
Here note that some or all of the entries $a_1b_1, a_1b_2, a_2b_1, a_2b_2$ may absent
in $D$. By the result on the numerical range of $2\times 2$ matrix,
we have the following proposition.

\begin{proposition} \label{3.2}
Let $a_1,a_2,b_1,b_2 \in \IC$ with $a_1 \ne a_2, b_1 \ne b_2$,
$c \in (0,1)$, $s = \sqrt{1-c^2}$, and
$$C = \begin{pmatrix}
a_1c^2+a_2s^2 & (a_1-a_2)cs \cr (a_1-a_2)cs & a_1 s^2 + a_2c^2\cr\end{pmatrix}
\begin{pmatrix} b_1 & 0 \cr 0 & b_2\cr\end{pmatrix}.$$
Then
\iffalse
$C$ has eigenvalues $[\gamma \pm \sqrt{\gamma^2 - 4a_1a_2b_1b_2}]/2$
and
\fi
$W(C)$ is the elliptical disk $\cE(a_1,a_2,b_1,b_2;\gamma)$
with foci $\gamma \pm \sqrt{\gamma^2 - a_1a_2b_1b_2}$
and length of minor axis
$$\{2|\hat\gamma|^2 + (|b_1|^2+|b_2|^2)|a_1-a_2|^2c^2s^2
- 2|\hat\gamma^2+b_1b_2(a_1-a_2)^2c^2s^2|\}^{1/2},$$
where $\gamma = \tr C = [(a_1b_1+a_2b_2)c^2+(a_1b_2+a_2b_1)s^2]$ and
 $\hat \gamma = [(a_1b_1-a_2b_2)c^2 + (a_2b_1-a_1b_2)s^2]/2$.
\end{proposition}

\medskip
Several remarks in connection to Proposition \ref{3.2} are in order.

\begin{enumerate}
\item
If $(a_1,a_2) = (b_1,b_2) = (1,0)$, then
$\cE(a_1,a_2,b_2,b_2;\gamma) = \cE(\gamma)$ defined in Section 2.
\item
The center of $W(C)$ in Proposition \ref{3.2} always lies in the line segment with
end points $a_1b_1 + a_2b_2$ and $a_1b_2+a_2b_1$,
and these two points are different if $a_1 \ne a_2$ and $b_1 \ne b_2$.

\item Suppose $a_1, a_2, b_1, b_2$ are given such that $a_1 \ne a_2, b_1 \ne b_2$.
Every $\gamma$ in the interior of the line segment with
end points $a_1b_1 + a_2b_2$ and $a_1b_2+a_2b_1$ uniquely determine $c \in (0,1)$ and
$s = \sqrt{1-c^2}$ so that one can construct the matrix $C$ (based on $a_1, a_2, b_1, b_2,
\gamma$) such that $W(C) = \cE(a_1,a_2,b_1,b_2;\gamma)$.

\item
Let $A, B \in M_n$ satisfy the hypothesis of Proposition \ref{3.1}.
Then for every
$\lambda \in (\sigma(AB) \setminus S)$
with
$S = \sigma(AB) \cap\{a_1b_1, a_1b_2, a_2b_1, a_2b_2\}$,
there is
$\tilde \lambda \in \sigma(AB)$ such that $\lambda\tilde \lambda = a_1a_2b_1b_2$
and $\lambda+\tilde \lambda = (a_1b_1+a_2b_2)c_j^2 + (a_1b_2+a_2b_1)s_j^2$.
Such a pair of eigenvalues correspond to the eigenvalues of $C_j$.
If $a_1a_2b_1b_2 = 0$, then
$\lambda \in (\sigma(C_j)\setminus S)$
will ensure that $\lambda \ne 0$ so that $\hat \lambda = 0$.
Otherwise, $\hat \lambda = a_1a_2b_1b_2/\lambda$. As a result,
we can always assume that $\hat \lambda = a_1a_2b_1b_2/\lambda$ and $\gamma = \lambda +
a_1a_2b_1b_2/\lambda$.
\end{enumerate}

By the above remarks and Propositions \ref{3.1}, \ref{3.2}, we have the following.

\begin{theorem} \label{3.3}
Suppose $A, B\in M_n$ are non-scalar normal matrices with
$\sigma(A) = \{a_1, a_2\}$ and $\sigma(B) = \{b_1, b_2\}$.
Let $S = \sigma(AB) \cap \{a_1b_1, a_1b_2, a_2b_1, a_2b_2\}$.
Then
$$W(AB) = \conv\{\cup_{\lambda \in (\sigma(AB)\setminus S)}
\cE(a_1,a_2,b_1,b_2;\lambda+ (a_1a_2b_1b_2)/\lambda) \cup S\}.$$
\end{theorem}

\medskip
We can use the dilation technique to study the numerical range of the product
of essentially Hermitian matrices.
Let $\tilde A$ be an essentially Hermitian matrix such that
$W(\tilde A)$ is a line segment joining $a_1, a_2 \in \IC$.
Then $\tilde A = a_2 I_n+ (a_1-a_2) A$ for a
positive contraction $A$.
Then $\tilde A$ has a dilation of the form
$\tilde P = a_2 I_{2n} + (a_1-a_2)P$
with
$$P = \begin{bmatrix} A & \sqrt{A-A^2} \cr \sqrt{A-A^2} & I-A\cr\end{bmatrix}.$$
Then $\tilde P$ is normal with $\sigma(\tilde P) = \{a_1, a_2\}$
so that $W(\tilde P) = W(A)$.
Now, if $\tilde B = b_2 I_n + (b_1-b_2)B$ is another essentially Hermitian
matrix, then $\tilde B$ has a dilation
$\tilde Q = b_2 I_{2n} + (b_1-b_2)Q$, where
$$Q = \begin{bmatrix} B & \sqrt{B-B^2} \cr \sqrt{B-B^2} & I-B\cr\end{bmatrix},$$
such that $\tilde Q$ is normal
with $\sigma(\tilde Q) = \{b_1, b_2\}$ and $W(\tilde B) = W(\tilde Q)$.

Using this observation and arguments similar to those in the proof of
Theorem \ref{2.2}, we have the following.

\begin{theorem} \label{3.4}
Suppose $A, B \in M_n$ are essentially Hermitian matrices such that
$A = a_2 I_n +  (a_1-a_2) A_1$
and $B = b_2 I_n + (b_1-b_2) B_1$ for two positive contractions $A_1, B_1$.
Let $\tilde A = a_2 I_{3n} + (a_1-a_2) P$ and
$\tilde B = b_2 I_{3n} + (b_1-b_2) Q$
$$P = \begin{bmatrix} A_1 & \sqrt{A_1-A_1^2} & 0 \cr \sqrt{A_1-A_1^2} & I_n-A_1 &0\cr
0 & 0 & 0\cr \end{bmatrix} \quad \hbox{ and }
\quad Q = \begin{bmatrix} B_1 & 0 & \sqrt{B_1 - B_1^2} \cr 0 & 0 & 0 \cr
\sqrt{B_1-B_1^2} & 0 & I_n-B_1 \cr \end{bmatrix}.$$
Then $W(AB) \subseteq W(\tilde A\tilde B)$, where $W(\tilde A\tilde B)$ can be determined by
Theorem \ref{3.3}. The set equality holds if and only if there is a unitary $U$ such that
$UAU^* = A_1 \oplus A_2, UBU^* = B_1 \oplus B_2$
satisfying $\sigma(A_1) = \{a_1, a_2\}$, $\sigma(B_1) = \{b_1, b_2\}$, and
$W(A_2B_2) \subseteq W(A_1B_1) = W(AB)$.
\end{theorem}

\medskip\noindent

\section{Extension to infinite dimensional spaces}

We can extend the results in the previous sections to $B(H)$, where $H$ is infinite dimensional.
Note that for a pair of non-scalar orthogonal projections
$P, Q \in B(H)$, there is a unitary $U$ such that
$U^*(P+iQ)U$ is a direct sum of $(1+i) I \oplus I \oplus iI \oplus 0$, and
$$\begin{bmatrix}C^2 + iI & C\sqrt{I-C^2} \cr C\sqrt{I-C^2} & I-C^2\cr\end{bmatrix},$$
where $C$ is a positive contraction; see \cite{Kl,S} and their references.
Consequently, $PQ$ is a direct sum of $I \oplus 0$ and
$$T = \begin{bmatrix}C^2 & 0 \cr C\sqrt{I-C^2} & 0\cr\end{bmatrix}.$$
Note that $T$
can be approximated by a sequence of operators of the form
$$T_m = \begin{bmatrix}C_m^2 & 0 \cr
C_m\sqrt{I-C_m^2} & 0 \cr\end{bmatrix}, \qquad m = 1, 2, \dots ,$$
where $C_m$ has finite spectrum and therefore can be assumed to be the direct
sum of $c_j^2 I_{H_j}$ on some subspace $H_j$ for $j = 1, \dots, c_{k_m}^2$
with $c_j \in (0,1)$. It is easy to see that
$$W(T_m) = \conv\{\cup_{\lambda \in \sigma(C_m)} \cE(\lambda)\}.$$
In fact, the same conclusion holds if $\sigma(C) = \sigma(T)$
is a discrete set, equivalently, $\sigma(PQ)$ is a discrete set.
In other words, if $P, Q \in B(H)$ are non-scalar orthogonal projections
such that $\sigma(PQ)$ is a finite or countably infinite, then
$$W(PQ) = \conv\{\cup_{\lambda \in \sigma(PQ)} \cE(\lambda)\}.$$
This result was proved in \cite[Theorem 1.3]{Kl},
for separable Hilbert spaces. One readily sees that the proof works for general
Hilbert space. In general, one can approximate $T$ by $T_m$ and
obtain the following result concerning the closure of $W(PQ)$; see \cite[Theorem 1.2]{Kl}.

\begin{proposition} \label{4.0}
Let $P, Q \in B(H)$ be non-scalar orthogonal projections. Then
$$\overline{W(PQ)} = \overline{\conv\{\cup_{\lambda \in \sigma(PQ)} \cE(\lambda)\}}.$$
The closure signs on both sides can be removed if $\sigma(PQ)$ is a discrete set.
\end{proposition}

One can show that the results
in Sections 2-3 hold in the infinite dimension setting.

\begin{theorem} \label{4.1}
Let $A,B \in B(H)$ be non-scalar positive semi-definite contractions.
Then
$$\overline{W(AB)} \subseteq \overline{\conv\{\cup_{\lambda \in \sigma(AB)} \cE(\lambda)\}}.$$
The closure signs can be removed if $\sigma(AB)$ is a discrete set.
The set equality holds if $A, B$ are orthogonal projections.
\end{theorem}

\begin{theorem} \label{4.2}
Let $A,B \in B(H)$ be such that $\overline{W(A)}$ is the line segment
joining $a_1, a_2$, and $\overline{W(B)}$ is the line segment joining $b_1, b_2$,
where $a_1 \ne a_2$ and $b_1 \ne b_2$.
Then
$$\overline{W(AB)} \subseteq
\overline{\conv\{\cup_{\lambda \in (\sigma(AB)\setminus S)}
\cE(a_1,a_2,b_1,b_2;\lambda+\hat \lambda) \cup S\}},$$
where $S= \sigma(AB) \cap \{a_1b_1, a_1b_2, a_2b_1, a_2b_2\}$,
$\hat \lambda = a_1a_2b_1b_2/\lambda$, and $\cE(a_1, a_2, b_1, b_2; \lambda+\hat \lambda)$
is defined as in Theorem \ref{3.3}.
The closure signs can be removed if $\sigma(AB)$ is a discrete set.
The set equality holds if $\sigma(A) = \{a_1, a_2\}$ and $\sigma(B) = \{b_1, b_2\}$.
\end{theorem}

\medskip
In Theorems \ref{4.1} and \ref{4.2}, we only have sufficient conditions for the set
inclusions become set equalities.
The problems of characterizing the set equality cases are open.

\bigskip\noindent
{\large\bf Acknowledgment}

The authors would like to thank Professor Ngai-Ching Wong for catalyzing the collaboration.
The research of Li was
supported by the USA NSF grant DMS 1331021, the Simons Foundation Grant 351047, and the
NNSF of China Grant 11571220.
The research of K.Z. Wang was supported by the Ministry of Science and Technology of the Republic
of China under project MOST 104-2918-I-009-001.
This research of Y. Wang and N. Zuo was partially supported by the Scholarship Program of
Chongqing University
of Science \& Technology, the Doctoral Research Fund of Chongqing University of Science
\& Technology (CK2010B09) and the National Natural Science Foundation of China (11571211).

\medskip
(H. Du) College of Mathematics and Information Science, Shaanxi Normal University,
Xi’an 710062, China. hkdu@snnu.edu.cn

(C.K. Li) Department of Mathematics, College of William and Mary,
Williamsburg, VA 23187, USA. ckli@math.wm.edu.

(K.Z. Wang) Department of Applied Mathematics,
National Chiao Tung University, Hsinchu 30010, Taiwan.
kzwang@math.nctu.edu.tw.

(Y. Wang) Department of Mathematics and Physics, Chongqing University of Science
and Technology, Chongqing 401331, China. wongyq@163.com

(N. Zuo) Department of Mathematics and Physics, Chongqing University of Science and
Technology, Chongqing 401331, China. zuon082@163.com

\end{document}